\documentclass[french]{article} 
\usepackage[utf8]{inputenc}
\usepackage[T1]{fontenc}
\usepackage{lmodern}
\usepackage[a4paper]{geometry}
\usepackage{fullpage}
\usepackage[english]{babel}

\usepackage{amssymb,amsmath,mathtools,amsthm,dsfont}

\usepackage[all]{xy}
\usepackage{ulem}

\newtheorem{theorem}{Theorem}
\newtheorem{prop-f}[theoreme]{Proposition}
\newtheorem{prop}[theorem]{Proposition}

\newtheorem{corollary}[theorem]{Corollary}

\newtheorem{lemma}[theorem]{Lemma}

\newcommand{\R}{\mathbb{R}}

\newcommand{\E}{\mathbb{E}}
\renewcommand{\P}{\mathbb{P}}
\renewcommand{\1}{\mathds{1}}

\renewcommand{\epsilon}{\varepsilon}
\renewcommand{\phi}{\varphi}

\newcounter{numeroexo}

\newcommand{\connectedto}{\overset{\Sigma}{\longleftrightarrow}}
\newcommand{\card}{\mbox{card}}

\title{\LARGE Equivalence of some subcritical properties\\ in continuum percolation}
\author{Jean-Baptiste Gou\'er\'e\footnote{Institut Denis-Poisson - UMR CNRS 7013, Universit\'e de Tours, Parc de Grandmont, 37200 Tours, France, {\it jean-baptiste.gouere@lmpt.univ-tours.fr}} \, and Marie Th\'eret\footnote{LPSM - UMR CNRS 8001, Universit\'e Paris Diderot, Sorbonne Paris Cit\'e, CNRS, F-75013 Paris, France, {\it marie.theret@univ-paris-diderot.fr}}}
\date{}

\begin{document}

\selectlanguage{\english}

\maketitle

\thispagestyle{empty}

\noindent
{\bf Abstract:} {We consider the Boolean model on $\R^d$. We prove some equivalences between subcritical percolation properties. Let us introduce some notations to state one of these equivalences. 
Let $C$ denote the connected component of the origin in the Boolean model. Let $|C|$ denotes its volume. Let $\ell$ denote the maximal length of a chain of random balls from the origin.  Under optimal integrability conditions on the radii, we prove that $E(|C|)$ is finite if and only if there exists $A,B >0$ such that $\P(\ell \ge n) \le Ae^{-Bn}$ for all $n \ge 1$.}
\\


\section{Introduction}	

\paragraph{The Boolean model.} The Boolean model is defined as follows.
At each point of a homogeneous Poisson point process on the Euclidean space $\R^d$, we center a ball of random radius.
We assume that the radii of the balls are independent, identically distributed and independent of the point process.
The Boolean model is the union of the balls.
There are three parameters: 
\begin{itemize}
\item An integer $d \ge 2$. This is the dimension of the ambient space $\R^d$.
\item A real number $\lambda>0$.  The intensity measure of the Poisson point process of centers is $\lambda |\cdot|$ where $|\cdot|$ denotes the Lebesgue measure on $\R^d$.
\item A probability measure $\nu$ on $(0,+\infty)$. This is the common distribution of the radii.
We will also consider a random variable $R$ whose distribution is $\nu$.
\end{itemize}
We will denote the Boolean model by $\Sigma(\lambda,\nu,d)$ or $\Sigma$.
We will also say that $\Sigma$ is the Boolean model driven by the measure $\lambda \nu$.

More precisely, the Boolean model is defined as follows. 
Let $\xi$ be a Poisson point process on $\R^d \times (0,+\infty)$ with intensity measure $\lambda|\cdot| \otimes \nu = |\cdot| \otimes \lambda \nu$.
Set
$$
\Sigma(\lambda,\nu,d)=\bigcup_{(c,r) \in \xi} B(c,r)
$$
where $B(c,r)$ denotes the open Euclidean ball of $\R^d$ with center $c$ and radius $r$.
We refer to the book by Meester and Roy \cite{Meester-Roy-livre} for background on the Boolean model, 
and to the book by Schneider and Weil \cite{Schneider-Weil} and the book by Last and Penrose \cite{last_penrose_2017} for background on Poisson processes.
We also denote by $S(c,r)$ the Euclidean sphere of $\R^d$ with center $c$ and radius $r$.
We write $S(r)$ when $c=0$.

\paragraph{Percolation in the Boolean model.} If $A$ and $B$ are two subsets of $\R^d$, we set
$$
\{A \connectedto B\} = \{\mbox{There exists a path in }\Sigma\mbox{ from }A\mbox{ to }B\}
$$
and
$$
\{0 \connectedto \infty\} = \{\mbox{The connected component of }\Sigma \mbox{ that contains the origin is unbounded}\}.
$$
We denote by $C$ the connected component of $\Sigma$ that contains the origin.
We denote by $D$ the diameter of $C$.
We denote by $\#C$ the number of random balls contained in $C$. 
In other words, $\#C$ is the following cardinality:
$$
\#C = \card(\{(c,r) \in \xi : c \in C\}).
$$
A chain of length $n \ge 1$ is a sequence $( (c_1,r_1), \dots, (c_n,r_n) )$ of distinct points of $\xi$ such that
$$
\forall i \in \{2, \dots, n\}, B(c_{i-1},r_{i-1}) \cap B(c_i,r_i) \neq\emptyset.
$$
We say that the chain starts in $A \subset \R^d$ if $B(c_1,r_1)$ touches $A$.
We say that the chain stops in $A \subset \R^d$ if $B(c_n,r_n)$ touches $A$.
We denote by $\ell$ the largest length of a chain starting in $B(0,1)$.
More precisely,
\begin{equation}\label{e:l}
\ell = \sup \{n \ge 0 : \exists x_1, \dots, x_n \in \xi \mbox{ s.t. } (x_1,\dots,x_n) \mbox{ is a chain starting in }B(0,1)\}.
\end{equation}

Let us consider the sets
\begin{align*}
 \Lambda  & = \{ \lambda > 0 : \P(0 \connectedto \infty)=0\}, \\
 \widehat{\Lambda} & = \{\lambda > 0 : \P(S(r) \connectedto S(2r)) \to 0 \mbox{ as } r \to \infty\}
\end{align*}
and the associated critical thresholds
\begin{align*}
 \lambda_c & = \sup \Lambda_c, \\
 \widehat{\lambda}_c & =  \sup \widehat{\Lambda}_c.
\end{align*}
From
$$
\P(0 \connectedto \infty) = \lim_{r \to \infty} \P(0 \connectedto S(2r))
$$
we get
$$
\widehat{\Lambda} \subset \Lambda
$$
and thus
$$
\widehat{\lambda}_c \le \lambda_c.
$$
When $\lambda \in \Lambda$, all the connected components of $\Sigma$ are almost surely bounded. We say that $\Sigma$ does not percolate.
When $\lambda \not\in \Lambda$, with probability one, one of the connected components of $\Sigma$ is unbounded. We say that $\Sigma$ percolates.
The thresholds $\lambda_c$ and $\widehat\lambda_c$ are always finite. 
See for example the remark below Theorem 3.3 in \cite{Meester-Roy-livre}.
If the radii are bounded, then $\lambda_c=\widehat\lambda_c$.
This is a sharp threshold property.
The sharpness of the transition in the discrete setting was proved independently 
by Menshikov \cite{Menshikov-coincidence} and by Aizenman-Barsky \cite{Aizenman-Barsky}.
The first proof of the equality $\lambda_c=\widehat\lambda_c$ relied on the analogous result in the discrete setting.
We refer to \cite{Meester-Roy-livre} for the proof (see Theorem 3.5) and references. Ziesche gives in \cite{Ziesche} a short proof of the equality $\lambda_c = \widehat \lambda_c$ for bounded radii. It relies on a new and short proof of the analogous result in the discrete setting by Duminil-Copin and Tassion \cite{DCT-sharp-Bernoulli, DCT-sharp-Ising}. Ahlbergh, Teixeira and Tassion gave in \cite{ATT} a very complete picture of percolation in the two dimensional Boolean model. In particular, they established a sharp threshold property for the two dimensional Boolean model under a minimal integrability assumption ($\E(R^2)<\infty$). See Theorem 1.1 in \cite{ATT}. We also refer to \cite{ATT2} by Ahlbergh, Teixeira and Tassion and \cite{Penrose17} by Penrose for further results about percolation in the complement of the Boolean model.
Even more recently, Duminil-Copin, Raoufi and Tassion developped new methods to prove sharp threshold properties in a wide class of models via decision trees, see for instance \cite{DC-Raoufi-Tassion-17-1,DC-Raoufi-Tassion-17-2}. In a course given at the IHES \cite{DC-video-17}, Duminil-Copin presented this new method and applied it to various models, including the Boolean model for which he annonces sharp transition in any dimension under some moment conditions.

Assume in this paragraph that $\E(R^d)$ is infinite.
Then, for any positive $\lambda$, with probability one, $\Sigma=\R^d$.
This can be shown easily by computing, for any $r>0$, the probability that $B(0,r)$ is covered by one random ball of $\Sigma$.
See for example Theorem 16.4 in \cite{last_penrose_2017}.
In this case, the model is therefore trivial from the percolation point of view: $\Lambda=\widehat{\Lambda}=\emptyset$ and $\lambda_c=\widehat\lambda_c=0$.
As a consequence, in what follows, we will always assume that $\E(R^d)$ is finite.

The following result is implicit in \cite{G-perco-boolean-model}.
A proof is given in the Appendix of \cite{Gouere-Theret-17} (Theorem 11).

\begin{theorem} \label{t:widehatpositif} Assume that $\E(R^d)$ is finite.
Then $\widehat \Lambda$ is open and non-empty.
In particular, $\widehat\lambda_c$ and $\lambda_c$ belong to $(0,+\infty)$.
\end{theorem}

In this paper, we want to investigate the connection between different percolation properties, such as the behavior of $\P(S(r) \connectedto S(2r))$ as $r $ goes to $\infty$, the integrability properties of $|C|$, $\#C$ and $D$, and the tail of the distribution of $\ell$. We state first the following result.

\begin{theorem} \label{t:implicit} Let $s>0$. The following statements are equivalent:
\begin{itemize}
 \item For small enough $\lambda$, $\E(|C|^{s/d})$ is finite.
 \item For small enough $\lambda$, $\E(\#C^{s/d})$ is finite.
 \item For small enough $\lambda$, $\E(D^s)$ is finite.
 \item $\E(R^{d+s})$ is finite.
\end{itemize}
Moreover, if $\E(R^{d+s})$ is finite, then $\E(|C|^{s/d})$, $\E(\#C^{s/d})$ and $\E(D^s)$ are finite as soon as $\lambda \in \widehat{\Lambda}$.
\end{theorem}

In particular, if $\E(R^{2d}) =\infty$ then $\E(|C|) = \E(\# C) = \E(D^d)=\infty$ for any $\lambda$. Since we are interested in the finiteness of those expectations, we will naturally suppose that $\E(R^{2d})<\infty$ in the following theorem, which is the main result of this article.

\begin{theorem} \label{t:new} 
Assume that $\E(R^{2d})$ is finite.
The following statements are equivalent:
\begin{enumerate}
 \item $\P(S(r) \connectedto S(2r)) \to 0$ as $r \to \infty$.
 \item There exists $A,B>0$ such that, for all $n \ge 1$,
 \begin{equation}\label{e:expo}
 \P( \ell \ge n ) \le A \exp(-Bn).
 \end{equation}
 \item $\E( D^d)$ is finite.
 \item $\E( |C| )$ is finite.
 \item $\E(\# C)$ is finite.
\end{enumerate}
\end{theorem}

The main contribution of our work is the proof of $4 \Rightarrow 2$, {\it i.e.}, the fact that $\E(|C|)<\infty$ implies the exponential decay of $\ell$, and this will be the core of the paper (see Section \ref{s:4-2}).
Note that \eqref{e:expo} does not imply that the decay of the tail of $\#C$ is exponential.
One can for example prove the following result, which is a simple consequence of Theorems \ref{t:implicit} and \ref{t:new}.

\begin{corollary} \label{c:1} 
Let $s>d$. Assume $\E(R^{2d})<\infty$ and $\E(R^{d+s})=\infty$.
Let $\lambda \in \widehat\Lambda = (0,\widehat{\lambda}_c)$. Then there exists $A,B$ such that
$$
\P( \ell \ge n ) \le A \exp(-Bn).
$$
However 
$$
\E(D^s)=\E(|C|^{s/d})=\E(\#C^{s/d}) = \infty.
$$
\end{corollary}

Combining Theorems \ref{t:implicit} and \ref{t:new} one also gets the following corollary.

\begin{corollary} \label{c:2} Let $s\ge d$. 
 Assume that $\E(R^{d+s})$ is finite.
The following statements are equivalent:
\begin{enumerate}
 \item $\P(S(r) \connectedto S(2r)) \to 0$ as $r \to \infty$.
 \item There exists $A,B>0$ such that, for all $n \ge 1$,
 $ \P( \ell \ge n ) \le A \exp(-Bn)$.
 \item $\E( D^s)$ is finite.
 \item $\E( |C|^{s/d} )$ is finite.
 \item $\E( \# C ^{s/d} )$ is finite.
\end{enumerate}
\end{corollary}

The above results also yield equalities between some percolation thresholds. Such equalities were already proven in the case where $R$ is bounded. We refer to Sections 3.4 and 3.5 of \cite{Meester-Roy-livre} and references therein for such results.

\paragraph{} The proof of Theorem \ref{t:implicit} is given in Section \ref{s:p:t:implicit}. The proof of Theorem  \ref{t:new} is given in Section \ref{s:p:t:new}. Corollaries \ref{c:1} and \ref{c:2} are straightforward consequences of Theorems \ref{t:implicit} and \ref{t:new}, thus no additional proof is needed.

\section{Proof of Theorem \ref{t:implicit}}
\label{s:p:t:implicit}

The proof of Theorem \ref{t:implicit} is divided into two parts. In Section \ref{s:p:t:implicit-1}, we prove the following result.
\begin{lemma}
\label{l:implicit-1}
Let $s>0$. If $\E[R^{d+s}] < \infty$ and $\lambda \in \widehat  \Lambda$, then $\E(|C|^{s/d})$, $\E(\#C^{s/d})$ and $\E(D^s)$ are finite. 
\end{lemma} 
Section \ref{s:p:t:implicit-2} is devoted to the proof of the following result.
\begin{lemma}
\label{l:implicit-2}
Let $s>0$. If $\E[R^{d+s}] = \infty$ then for every $\lambda>0$, $\E(|C|^{s/d})$, $\E(\#C^{s/d})$ and $\E(D^s)$ are infinite. 
\end{lemma}
Theorem \ref{t:implicit} is a straightforward consequence of Lemmas \ref{l:implicit-1} and \ref{l:implicit-2}.


\subsection{Proof of Lemma \ref{l:implicit-1}}
\label{s:p:t:implicit-1}

We first establish the following result.
\begin{theorem} \label{t:explicitons} Let $s>0$. Assume $\E(R^{d+s})<\infty$.
Let $\lambda \in \widehat{\Lambda}$.
Then
$$
\int_0^\infty \alpha^{s-1} \P(S(\alpha) \connectedto S(2\alpha)) \; d\alpha < \infty
$$
and
$$
\E(D^s)<\infty.
$$
\end{theorem}

The result is implicit in \cite{G-perco-boolean-model}.
We choose to give a detailed proof using intermediate results in Appendix A in \cite{Gouere-Theret-17} which themselves rely on results in \cite{G-perco-boolean-model}.

\bigskip

Let us recall some notation from \cite{G-perco-boolean-model} or \cite{Gouere-Theret-17}.
Let $\alpha>0$.  
\begin{itemize}
 \item $\Sigma(B(0,\alpha))$ is the union of random balls of the Boolean model with centers in $B(0,\alpha)$.
 \item $G(0,\alpha)$ is the event "there exists a path from $S(\alpha)$ to $S(8\alpha)$ in $\Sigma(B(0,10\alpha))$".
 \item $\Pi(\alpha)=P(G(0,\alpha))$.
\end{itemize}
Set 
\begin{equation} \label{e:application-epsilon}
 \epsilon(\alpha)=\int_{[\alpha,+\infty)} r^d \nu(dr).
\end{equation}
Note that, when $\E(R^{d+s})$ is finite, 
\begin{equation}\label{e:capiquetoujours}
\int_0^\infty \alpha^{s-1}\epsilon(\alpha) \;d\alpha < \infty.
\end{equation}

The following proposition is stated in the same way in \cite{Gouere-Theret-17} as Proposition 12 in Appendix A.

\begin{prop} \label{l:prop} There exists a constant $K=K(d)$ such that, for any $\alpha > 0$,
\begin{eqnarray} 
\Pi(\alpha) & \le & P(S(\alpha) \connectedto S(2\alpha)) \le K\Pi(\alpha/10) + \lambda K \epsilon(\alpha/10), \label{e:olivier} \\
\Pi(10\alpha) & \le & K\Pi(\alpha)^2+\lambda K \epsilon(\alpha), \label{e:florence} \\
\Pi(\alpha) & \le & \lambda K \alpha^d \label{e:claire}.
 \end{eqnarray}
\end{prop}

\paragraph{Proof of Theorem \ref{t:explicitons}.} This is a consequence of Proposition \ref{l:prop} above and Lemma 3.7 in \cite{G-perco-boolean-model}. 
Showing how to apply Lemma 3.7 would not be much shorter than adapting the proof in our context. 
Therefore we choose to give a full proof.
Let $s>0$. Assume $\E(R^{d+s})<\infty$ and let $\lambda \in \widehat\Lambda$.
By \eqref{e:olivier}, $\Pi(\alpha)$ tends to $0$ as $\alpha$ tends to $\infty$ 
Therefore we can fix $\alpha_0$ large enough such that, for all $\alpha \ge \alpha_0/10$,
$$
10^sK\Pi(\alpha) \le \frac 12.
$$
Then, for any $\alpha \ge \alpha_0$, using \eqref{e:florence} and the definition of $\alpha_0$,
\begin{align*}
 \int_{\alpha_0}^\alpha r^{s-1}\Pi(r)dr 
  & \le \int_{\alpha_0}^\alpha r^{s-1} K\Pi(r/10)^2dr + \int_{\alpha_0}^\alpha r^{s-1}\lambda K \epsilon(r/10) dr \\
  & \le 10^s\int_{\alpha_0/10}^{\alpha/10} r^{s-1} K\Pi(r)^2dr + \int_{\alpha_0}^\infty r^{s-1}\lambda K \epsilon(r/10) dr \\
  & \le \frac 12 \int_{\alpha_0/10}^{\alpha/10} r^{s-1} \Pi(r)dr + \int_{\alpha_0}^\infty r^{s-1}\lambda K \epsilon(r/10) dr.
\end{align*}
Therefore, for any large enough $\alpha$,
$$
 \int_{\alpha_0}^\alpha r^{s-1}\Pi(r)dr \le \frac 12 \int_{\alpha_0/10}^{\alpha_0} r^{s-1} \Pi(r)dr + \frac 12 \int_{\alpha_0}^\alpha r^{s-1} \Pi(r)dr + \int_{\alpha_0}^\infty r^{s-1}\lambda K \epsilon(r/10) dr.
$$
Then, rearranging and using \eqref{e:capiquetoujours},
$$
\int_{\alpha_0}^\alpha r^{s-1}\Pi(r)dr \le \int_{\alpha_0/10}^{\alpha_0} r^{s-1} \Pi(r)dr + 2\int_{\alpha_0}^\infty r^{s-1}\lambda K \epsilon(r/10) dr < \infty.
$$
Therefore,
$$
\int_{\alpha_0}^\infty r^{s-1}\Pi(r)dr < \infty
$$
and thus
$$
\int_0^\infty r^{s-1}\Pi(r)dr < \infty.
$$
By \eqref{e:capiquetoujours} and \eqref{e:olivier}, this yields the first required result.
The other result then follows from the fact that, for any $\alpha>0$,
$$
\{D \ge 4\alpha\} \subset \{S(\alpha) \connectedto S(2\alpha)\}.
$$

\qed

\paragraph{Proof of Lemma \ref{l:implicit-1}.}
We suppose that $\E(R^{d+s}) <\infty$ and $\lambda \in \widehat{\Lambda}$. By Theorem \ref{t:explicitons}, we know that $\E(D^s)<\infty$. Since $C\subset B (0, D)$, this implies $\E(|C|^{s/d})<\infty$. It remains to prove that $\E(\#C^{s/d})< \infty$.

Let $\kappa>0$ be such that
$$
\lambda v_d \kappa^d = \frac 1 2,
$$
where $v_d$ denotes the volume of the unit ball in $\R^d$.
For every $u>0$, we have
\begin{align}
\label{e:ajout}
\P(\#C \ge u ) 
 & \le  \P( C \subset B(0,\kappa u^{1/d}) \text{ and } \#C \ge u) + \P( C \not\subset B(0,\kappa u^{1/d})) \nonumber \\
 & \le  \P( \#\{(c,r) \in \xi : c \in B(0,\kappa u^{1/d})\} \ge u) + \P( D \ge \kappa u^{1/d}) .
\end{align}
Since $\#\{(c,r) \in \xi : c \in B(0,\kappa u^{1/d})\}$ is a Poisson random variable with parameter $u/2$, we obtain
$$
\P( \#\{(c,r) \in \xi : c \in B(0,\kappa u^{1/d})\} \ge u) \le \exp\left(u\left(\frac 1 2 - \ln(2) \right)\right)
$$
and thus
$$
\int_0^\infty du \; \frac s d u^{\frac s d - 1} \P( \#\{(c,r) \in \xi : c \in B(0,\kappa u^{1/d})\} \ge u) < \infty.
$$
Since $\E(D^s)<\infty$, we have
$$
\int_0^\infty du \; \frac s d u^{\frac s d - 1} \P(D \ge \kappa u^{1/d}) < \infty.
$$
We conclude by \eqref{e:ajout} that $\E(\# C^{s/d})$ is finite. \qed


\subsection{Proof of Lemma \ref{l:implicit-2}}
\label{s:p:t:implicit-2}

Set 
$$
A = \sup \{ r > 0 : \exists c \in B(0,r/2) \mbox{ s.t.} (c,r) \in \xi\},
$$
with the convention $A=0$ if the set is empty. Note that $B(0,A/2)$ is covered by $\Sigma$. We first state the following preliminary result, which is essentially implicit in \cite{G-perco-boolean-model}.

\begin{lemma}\label{l:A} Let $\lambda>0$ and $s>0$. Assume that $\E(R^{d+s})$ is infinite. Then $\E(A^s)$ is infinite.
\end{lemma}

\paragraph{Proof of Lemma \ref{l:A}.} For any $a>0$,
$$
\P(A > a) = 1-\exp\left(-\lambda v_d 2^{-d} \int_{(a,+\infty)} r^d \nu(dr)\right).
$$
If $\E(R^d)$ is infinite, then $\P(A>a)=1$ for all $a>0$ and thus $A=+\infty$ almost surely and therefore $\E(A^s)=\infty$.

Assume henceforth that $\E(R^d)$ is finite. Then
$$
\P(A>a) \sim_{a \to \infty} \lambda v_d 2^{-d} \int_{(a,+\infty)} r^d \nu(dr).
$$
Therefore, for some constant $\gamma>0$, 
for all $a>0$,
$$
\P(A > a) \ge \gamma \int_{(a,+\infty)} r^d \nu(dr)
$$
and then
\begin{align*}
\E(A^s) 
 & =  \int_{(0,+\infty)} da\;  sa^{s-1} \P(A >a) \\ 
 & \ge \gamma   \int_{(0,+\infty)} da \;sa^{s -1} \int_{(a,+\infty)} \nu(dr) \; r^d \\
 & = \gamma   \int_{(0,+\infty)} \nu(dr) \; r^d \int_{(0,r)} da \;s a^{s-1} \\
 & = \gamma   \int_{(0,+\infty)} \nu(dr) \; r^{d+s} \\
 & = \gamma \E(R^{s+d})
\end{align*}
which is infinite by assumption. \qed

\paragraph{Proof of Lemma \ref{l:implicit-2}.}
Let $\lambda>0$ and $s>0$. We suppose that $\E(R^{d+s})$ is infinite. By Lemma \ref{l:A}, we obtain that $\E(A^s)=\infty$. Since $B(0,A/2)$ is covered by $\Sigma$, we know that $D\geq A$ and $|C|\geq v_d A^d / 2^d$, thus $\E (D^s) = \infty$ and $\E(|C|^{s/d}) = \infty$. It remains to prove that $\E(\# C ^{s/d}) = \infty$.

Let $r_0>0$ be such that  $\P(R \le r_0)>0$.
Set
$$
A_{>r_0} = A \mbox{ if } A>r_0 \mbox{ and } A_{>r_0} =0 \mbox{ otherwise.}
$$
In other words,
$$
A_{>r_0} = \sup \{ r > r_0 : \exists c \in B(0,r/2) \mbox{ s.t. } (c,r) \in \xi\}
$$
with the convention $A_{>r_0}=0$ if the set is empty.
Note that $B(0,A_{>r_0}/2)$ is covered by $\Sigma$ and that $A_{>r_0}$ is measurable with respect to  $\xi_{> r_0} = \xi \cap \R^d  \times (r_0,+\infty)$.
Set 
$$
N = \card( \{ (c,r) \in \xi : c \in B(0,A_{>r_0}/2) \mbox{ and } r \le r_0 \}).
$$
Conditionally on $\xi_{> r_0}$ 
, $N$ is a Poisson random variable with parameter $\alpha A_{>r_0}^d$ where 
$$
\alpha=\lambda \P(R \le r_0) v_d / 2^d >0.
$$
But as $B(0,A_{>r_0}/2)$ is contained in $\Sigma$, any random ball centered in $B(0,A_{>r_0}/2)$ is contained in $C$.
Therefore
$$
\#C \ge N
$$
and thus 
$$
\E(\#C^{s/d}) \ge \E(N^{s/d}) = \E(\E(N^{s/d} | \xi_{> r_0})).
$$
Let $\mu_0$ be such that, for any $\mu \ge \mu_0$, if $X(\mu)$ is a Poisson random variable with parameter $\mu$, 
then $\P(X(\mu) \ge \mu/2) \ge 1/2$.
Then, 
\begin{align*}
\E(\#C^{s/d}) & \ge  \E\left(\E(N^{s/d} | \xi_{> r_0}) \1_{\alpha A_{>r_0}^d \ge \mu_0}\right) \\
& \ge \E\left(\frac 1 2 \left(\frac{\alpha A_{>r_0}^d}2\right)^{s/d}  \1_{\alpha A_{>r_0}^d \ge \mu_0}\right)
 = \beta \E(A_{>r_0}^s \1_{\alpha A_{>r_0}^d \ge \mu_0})
\end{align*}
where
$$
\beta = \frac{\alpha^{s/d}}{2^{1+s/d}}.
$$
To end the proof, it is therefore sufficient to check that $\E(A_{>r_0}^s)$ is infinite. 
But this is a consequence of the infiniteness of $\E(A^s)$ obtained by Lemma \ref{l:A}.
\qed


\section{Proof of Theorem  \ref{t:new}}
\label{s:p:t:new}

\subsection{Preliminary results.} 

For any $r>0$, we denote by $C(r)$ the connected component of $\Sigma \cup B(0,r)$ which contains $B(0,r)$.
If $A$ and $B$ are two subsets of $\R^d$, we denote by $A+B$ the Minkowski sum of $A$ and $B$ defined by
$$
A+B=\{a+b, a \in A, b \in B\}.
$$

\begin{lemma} \label{l:finalementfacile} Let $\lambda>0$. Assume that $\E(|C|)$ is finite. Then
 $$
  \E(|C(1)+B(0,1)|)<\infty.
 $$
\end{lemma}

\paragraph{Proof.} 
This is a consequence of FKG inequality, see for instance Theorem 2.2 in \cite{Meester-Roy-livre}.
Set 
$$
p = \P( B(0,1) \subset \Sigma) 
$$
and note that $p$ is positive.
For any $x \in \R^d$, using FKG inequality in the third step and stationarity and definition of $p$ in the fourth step, we get
\begin{align*}
 \P(x \in C) & = \P(\{0\} \connectedto \{x\}) \\
& \ge \P(\{B(0,1) \connectedto B(x,1) \} \cap \{B(0,1) \subset \Sigma\} \cap \{B(x,1) \subset \Sigma\}\big) \\
& \ge \P(B(0,1) \connectedto B(x,1))\P(B(0,1) \subset \Sigma)\P(B(x,1) \subset \Sigma) \\
& = p^2\P(B(0,1) \connectedto B(x,1)) \\
& = p^2\P(x \in C(1)+B(0,1)).
\end{align*}
Therefore 
$$
\E(|C(1)+B(0,1)|) = \int_{\R^d} dx \; \P(x \in C(1)+B(0,1)) \le p^{-2} \int_{\R^d} dx \; \P(x \in C) = p^{-2} \E(|C|).
$$
As $\E(|C|)$ is finite, the lemma is proven. \qed

\begin{lemma} \label{l:geo} 
Let $K=K(d)$ be such that, for any $r>0$, the ball $B(0,r)$ can be covered by $K(1+r)^d$ balls of radius $1$. 
Let $\lambda>0$. Let $r,s>0$. Then 
$$
\int_{\R^d} dx \; \P(B(0,r) \connectedto B(x,s)) = \E(|C(r)+B(0,s)|) \le K^2(1+r)^d(1+s)^d \E(|C(1)+B(0,1)|).
$$
\end{lemma}

\paragraph{Proof.} Write
$$
\E(|C(r)+B(0,s)|) = \int_{\R^d} dx \; \P(x \in C(r) + B(0,s)) = \int_{\R^d} dx \; \P(B(0,r) \connectedto B(x,s)).
$$
Cover $B(0,r)$ with at most $K(1+r)^d$ balls of radius $1$.
Cover $B(x,s)$ with at most $K(1+s)^d$ balls of radius $1$.
If $\{\P(B(0,r) \connectedto B(x,s)\}$ holds, then there exists a path in $\Sigma$ from one of the balls that cover $B(0,r)$ to one of the balls that cover $B(x,s)$.
Using union bound, stationarity and a change of variable, we thus get
\begin{align*}
\E(|C(r)+B(0,s)|) 
 & \le K^2(1+r)^2(1+s)^2 \int_{\R^d} dx \; \P(B(0,1) \connectedto B(x,1)) \\
 & = K^2(1+r)^d(1+s)^d \E(|C(1)+B(0,1)|). 
\end{align*}
The lemma is proven. \qed

\subsection{Proof of $4 \implies 2$} 
\label{s:4-2}

This is the main part of the proof. For any $n \ge 2$, any $r,s>0$ and any $x,y \in \R^d$ we consider the event
\begin{align*}
 L_n(x,r,y,s) &  = \{\mbox{There exists a chain of length } n-1 \mbox{ starting in } B(x,r) \mbox{ and stopping in } B(y,s)\} \\
\end{align*}
and we set
\begin{equation} \label{e:a}
a_n(r,s)= \int_{\R^d} dy \;  \P(L_n(0,r,y,s)).
\end{equation}
Recall that $R$ is a random variable with distribution $\nu$.
Let $S_i, i \ge 0$ be independent copies of $R+1$.

\begin{lemma} \label{l:debut} For any $\lambda>0$, any $p \ge 2$ and any $k \ge 1$,
$$
\P(\ell \ge kp)  \le \left(\lambda \left[\E\left(a_p(S_1,S_2)^2\right)\right]^{1/2}\right)^k.
$$
\end{lemma}

\paragraph{Proof.}
Set $c_0=0$ and $r_0=1$.
Let $\lambda>0$, $p \ge 2$ and $k \ge 1$.
We will use BK inequality, see the main theorem in \cite{vdB-BK} (and the remark (iii) above it concerning the choice of the definition of disjoint occurence of increasing events).
Before stating the inequality, let us recall informally some notations. The BK inequality apply directly to increasing events living on bounded region, thus define
\begin{align*}
 L^n_p(x,r,y,s) &  = \left\{ \begin{array}{c} \mbox{There exists a chain } ((x_1,r_1),\dots , (x_{p-1}, r_{p-1})) \mbox{ of length } p-1 \\ \mbox{ starting in } B(x,r)  \mbox{ and stopping in } B(y,s) \\\mbox{ s.t. for all } i\in \{ 1,\dots, p-1 \}, \, x_i \in [-n,n]^d\end{array} \right\}. \\
\end{align*}
If $((c_i,r_i))_{1\leq i \leq k}$ are points of $\R^d \times (0,+\infty)$, we say that the increasing events $L_p(c_{i-1},r_{i-1},c_i,r_i)$ (respectively $L^n_p(c_{i-1},r_{i-1},c_i,r_i)$), $i \in \{1,\dots,k\}$, occur disjointly
if there exists $k$ chains, each of length $p-1$, using in total $k(p-1)$ distinct random balls 
such that, for all $i\in\{1,\dots,k\}$ the $i$-th chain starts in $B(c_{i-1},r_{i-1})$ and stops in $B(x_i,r_i)$ (respectively and all the centers of the balls of these chains belong to $[-n,n]^d$).
We denote these events by
\begin{equation*}
L_p(c_0,r_0,c_1,r_1) \circ \cdots \circ L_p(c_{k-1},r_{k-1},c_k,r_k) \mbox{ and } L^n_p(c_0,r_0,c_1,r_1) \circ \cdots \circ L^n_p(c_{k-1},r_{k-1},c_k,r_k)
\end{equation*}
or simply by
\begin{equation*}
\underset{i}{\circ} \; L_p(c_{i-1},r_{i-1},c_i,r_i) \mbox{ and } \underset{i}{\circ} \; L^n_p(c_{i-1},r_{i-1},c_i,r_i).
\end{equation*}
By BK inequality, for all $n\in \mathbb N$ we have
\begin{equation*}
\P(L^n_p(c_0,r_0,c_1,r_1) \circ \cdots \circ L^n_p(c_{k-1},r_{k-1},c_k,r_k)) \le \prod_{i=1}^k \P(L^n_p(c_{i-1},r_{i-1},c_i,r_i	)).
\end{equation*}
Taking the limit as $n$ goes to infinity, we obtain
\begin{equation} \label{e:BK}
\P(L_p(c_0,r_0,c_1,r_1) \circ \cdots \circ L_p(c_{k-1},r_{k-1},c_k,r_k)) \le \prod_{i=1}^k \P(L_p(c_{i-1},r_{i-1},c_i,r_i	)).
\end{equation}
If $\ell \ge kp$, then there exists a chain of $kp$ distinct balls starting in $B(0,1)$.
Taking one ball every $p$-th balls in this chain, we get a sequence $(c_1,r_1),\dots,(c_k,r_k)$ of distinct points of $\xi$ such that, with a slight abuse of notation, the event
$$
\underset{i}{\circ} \; L_p(c_{i-1},r_{i-1},c_i,r_i)
$$
holds for $\xi \setminus \{(c_1,r_1),\dots,(c_k,r_k)\}$.
Therefore, again with a slight abuse of notation,
$$
\P(\ell \ge kp)  \le \E\left(\sum_{(c_1,r_1),\dots,(c_k,r_k) \in \xi \mbox{ distinct}} \1_{\underset{i}{\circ} \; L_p(c_{i-1},r_{i-1},c_i,r_i)}(\xi \setminus \{(c_1,r_1),\dots,(c_k,r_k)\}) \right).
$$
From Slivnyak’s theorem, see Proposition 4.1.1 in \cite{Moller}, we get
$$
\P(\ell \ge kp)  \le  \lambda^k \int_{(\R^d)^k} dc_1\dots dc_k \int_{(0,+\infty)^d} \nu(dr_1)\dots\nu(dr_k) \; \P\left(\underset{i}{\circ} \; L_p(c_{i-1},r_{i-1},c_i,r_i)\right).
$$
By \eqref{e:BK}, this yields
$$
\P(\ell \ge kp)  \le  \lambda^k  \int_{(\R^d)^k} dc_1\dots dc_k \int_{(0,+\infty)^d} \nu(dr_1)\dots\nu(dr_k) \prod_{i=1}^k \P(L_p(c_{i-1},r_{i-1},c_i,r_i)).
$$
Using stationarity and \eqref{e:a} we then get
$$
\P(\ell \ge kp)  \le  \lambda^k  \int_{(0,+\infty)^d} \nu(dr_1)\dots\nu(dr_k) \prod_{i=1}^k a_p(r_{i-1},r_i).
$$
Distinguishing according to parity, we get
$$
\P(\ell \ge kp)  \le  \lambda^k  \int_{(0,+\infty)^d} \nu(dr_1)\dots\nu(dr_d) \prod_{1 \le i \le k, i \mbox{ odd}} a_p(r_{i-1},r_i)  \prod_{1 \le i \le k, i \mbox{ even}} a_p(r_{i-1},r_i).
$$
Then, by Cauchy-Schwarz inequality, 
$$
\P(L_{kp}(1))  \le  \lambda^k  A_p(k)B_p(k)
$$
where
$$
A_p(k)=\left[\int_{(0,+\infty)^d} \nu(dr_1)\dots\nu(dr_d) \prod_{1 \le i \le k, \; i \mbox{ odd}} a_p(r_{i-1},r_i)^2\right]^{1/2}
$$
and
$$
B_p(k)=\left[\int_{(0,+\infty)^d} \nu(dr_1)\dots\nu(dr_d) \prod_{1 \le i \le k, \; i \mbox{ even}} a_p(r_{i-1},r_i)^2\right]^{1/2}.
$$
Let us get rid of the (easy but annoying) special case $r_0=1$ as follows.
Recall that $R$ is a random variable with distribution $\nu$ and that $S_i, i \ge 0$ are independent copies of $R+1$.
As $a_p(r,s)$ is non-decreasing in $r$ and $s$, we have
$$
A_p(k) \le \left[\E\left(\prod_{1 \le i \le k, \; i \mbox{ odd}} a_p(S_{i-1},S_i)^2\right)\right]^{1/2}
$$
and
$$
B_p(k) \le \left[\E\left(\prod_{1 \le i \le k, \; i \mbox{ even}} a_p(S_{i-1},S_i)^2\right)\right]^{1/2}.
$$
As $(S_i)_{i \ge 0}$ is an i.i.d. sequence, we get
$$
A_p(k)B_p(k) \le \left[\E\left(a_p(S_1,S_2)^2\right)\right]^{k/2}
$$
and then
$$
P(\ell \ge kp) \le \left(\lambda \left[\E\left(a_p(S_1,S_2)^2\right)\right]^{1/2}\right)^k.
$$
The lemma is proven. \qed

\paragraph{} Recall that $\ell$ is defined by $\eqref{e:l}$.

\begin{lemma} \label{l:trou} Let $\lambda>0$. Assume that $\E(R^d)$ and $\E(|C|)$ are finite. Then $\ell$ is finite with probability one.
\end{lemma}

\paragraph{Remark. } We could remove the assumption $\E(R^d)$ finite, as it is a consequence of the finiteness of $\E(|C|)$.

\paragraph{Proof.} Let $\lambda>0$ be such that $\ell$ is infinite with positive probability. We aim at proving that $\E(|C|)$ is infinite.
For any $\eta>0$, set
$$
\Sigma_{\le \eta} = \bigcup_{(c,r) \in \xi : r \le \eta} B(c,r).
$$
If $\eta$ is small enough, then $\Sigma_{\le \eta}$ does not percolate.
Indeed, $\Sigma_{\le \eta}$ is a Boolean model driven by the measure
$$
\lambda \nu ( \cdot \cap (0,\eta] ).
$$
We can therefore couple $\Sigma_{\le \eta}$ with the Boolean model $\Sigma^+_{\le \eta}$ driven by the measure
$$
\lambda \nu((0,\eta]) \delta_\eta
$$ 
in such a way that
$$
\Sigma_{\le \eta} \subset \Sigma^+_{\le \eta}.
$$
But 
$$
\eta^{-1} \Sigma^+_{\le \eta}
$$
is a Boolean model driven by 	
$$
\lambda \nu((0,\eta]) \eta^d \delta_1.
$$
Therefore, as soon as
$$
\lambda \nu((0,\eta]) \eta^d < \lambda_c(\delta_1,d),
$$
$\eta^{-1} \Sigma^+_{\le \eta}$ does not percolate, 
thus $\Sigma^+_{\le \eta}$ does not percolate and then $\Sigma_{\le \eta}$ does not percolate.
In the remaining of the proof, we fix $\eta>0$  such that $\Sigma_{\le \eta}$ does not percolate.

Let $C_{>\eta}$ denote the union of all the random balls of $C$ whose radius is greater than $\eta$. We will repeatedly use the following property, that holds since $\E(R^d)$ is finite: almost surely, whatever the bounded region $B$ of $\mathbb R^d$ we consider, the number of random balls of $\Sigma$ that touches $B$ is finite. If $\ell$ is infinite, then for any $n$ there exists a chain $((c^n_1, r^n_1), \dots , (c^n_n,r^n_n))$ of balls starting in $B(0,1)$. For all $n$ the ball $B(c^n_1, r^n_1)$ touches $B(0,1)$, but by the previously stated property we know that a.s. there exists only a finite number of balls of $\Sigma$ that touches $B(0,1)$, thus an infinite number of those balls $B(c^n_1, r^n_1)$ are equal to the same ball that we will denote by $B_1 $. In other words, there exists an infinite number of chains of arbitrarily large length starting in $B_1$. But by the same property, we know that a.s. there exists only a finite number of balls of $\Sigma$ that touches $B(0,1) \cup B_1$, thus an infinite number of those chains (with arbitrarily large length) use a common first ball of $\Sigma$ that we denote by $B_2$. By induction, we construct an infinite chain of distinct balls $(B_i = B(x_i, s_i) , i\geq 1)$ starting in $B(0,1)$. Suppose that only a finite number of the radii $s_i$ are bigger than $\eta$. Then there exists an infinite chain of balls of radii smaller than or equal to $\eta$. By the previously stated property, we know that a.s. this infinite number of balls cannot stay in any bounded region, thus $\Sigma_{\leq \eta}$ has to percolate, which is absurd by our choice of $\eta$. We conclude that an infinite number of these balls $(B_i = B(x_i, s_i) , i\geq 1)$ satisfy $s_i> \eta$. A.s. these balls cannot stay in any bounded region, thus up to extraction we obtain a sequence $(c_i,r_i)_{i \ge 0}$ of points of $\xi$ such that, for all $i \ge 1$, 
$$
B(c_i,r_i) \subset C, \quad \|c_i\| \ge \|c_{i-1}\|+1, \quad r_i \ge \eta.
$$
Therefore $|C|=\infty$ almost surely on the event $\{\ell = \infty\}$.
As this event occurs with positive probability we get $\E(|C|)=\infty$. \qed

\begin{lemma} \label{l:fin} Let $\lambda>0$. Assume that $\E(R^{2d})$ and $\E(|C|)$ are finite. Then there exists $p \ge 2$ such that
$$
\lambda \left[\E\left(a_p(S_1,S_2)^2\right)\right]^{1/2} < 1.
$$
\end{lemma}

\paragraph{Proof.} As $\E(|C|)$ is finite, $\ell$ is almost surely finite by Lemma \ref{l:trou}. 
For any $r>0$, $B(0,r)$ can be covered by a finite number of balls of radius $1$.
Therefore, by stationarity, the maximal number $\ell(r)$ of balls in a chain starting in $B(0,r)$ is almost surely finite.
As a consequence, for all $r>0$,
$$
\P(\ell(r) \ge p) \to 0 \mbox{ as } p \to \infty.
$$
Thus, for all $r,s>0$ and $y \in \R^d$,
$$
\P(L_p(0,r,y,s)) \to 0 \mbox{ as } p \to \infty,
$$
as $L_p(0,r,y,s) \subset \{\ell(r) \ge p-1\}$.
Moreover, for all  $r,s>0$, $y \in \R^d$ and $p \ge 2$,
\begin{equation} \label{e:13h}
\P(L_p(0,r,y,s)) \le \P(B(0,r) \connectedto B(y,s)).
\end{equation}
But, by Lemma \ref{l:geo},
\begin{equation}\label{e:reunionmaquette}
\int_{\R^d} dy \; \P(B(0,r) \connectedto B(y,s)) \le K^2(1+r)^d(1+s)^d \E(|C(1)+B(0,1)|).
\end{equation}
Moreover, as $\E(|C|)$ is finite, we get
\begin{equation}\label{e:maintenant}
\E(|C(1)+B(0,1)|) < \infty
\end{equation}
by Lemma \ref{l:finalementfacile}.
Therefore, by dominated convergence theorem, for any $r,s>0$,
$$
\int_{\R^d} dy \; \P(L_p(0,r,y,s))  \to 0 \mbox{ as } p \to \infty
$$
that is, using Definition \eqref{e:a},
$$
a_p(r,s)  \to 0 \mbox{ as } p \to \infty.
$$
By \eqref{e:13h} and \eqref{e:reunionmaquette},
$$
a_p(S_1,S_2) \le K^2(1+S_1)^d(1+S_2)^d \E(|C(1)+B(0,1)|).
$$
Using \eqref{e:maintenant} and the finiteness of $\E(R^{2d})$ we get that the square of the right hand side of the above inequality is integrable.
Using dominated convergence theorem again, we then get
$$
\E(a^2_p(S_1,S_2)) \to 0 \mbox{ as } p \to \infty.
$$
The lemma is proven. \qed

\paragraph{Proof of $4 \implies 2$.} By Lemma \ref{l:fin}, there exits $p \ge 2$ such that
$$
\kappa := \lambda \left[\E\left(a_p(S_1,S_2)^2\right)\right]^{1/2} < 1.
$$
By Lemma \ref{l:debut}, for any $k \ge 1$,
$$
\P(\ell \ge kp) \le \kappa^k.
$$
Therefore there exists $A,B>0$ such that, for any $n \ge 1$, $\P(\ell \ge n) \le A \exp(-Bn)$. \qed

\subsection{Proof of the others implications}

\paragraph{Proof of $2 \implies 1$.} Let $K=K(d)$ be such that, for any $r \ge 1$, the sphere $S(r)$ can be covered by $Kr^{d-1}$ balls of radius $1$.
Let $r \ge 1$. 
Let us first prove
\begin{equation}\label{e:poildechameau}
\P(S(r) \connectedto S(2r)) \le \P(H(r)) + Kr^{d-1}\P(\ell \ge \sqrt{r}/2)
\end{equation}
where
$$
H(r) = \{ \exists (c,s) \in \xi : s \ge \sqrt{r} \mbox{ and } B(c,s) \cap B(0,2r) \neq \emptyset\}.
$$
Indeed, if $\{S(r) \connectedto S(2r)\}$ holds, then there exists a chain of random balls from a point of $S(r)$ to a point of $S(2r)$.
If moreover $H(r)$ does not hold, then the number of balls of this chain is at least $\sqrt{r}/2$.
The sphere $S(r)$ can be covered by $Kr^{d-1}$ balls of radius $1$ and the starting point of the chain is in one of these balls.
Using the union bound and stationarity, we get \eqref{e:poildechameau}.

As
$$
\P(\ell \ge \sqrt{r}/2) \le A \exp(-B \sqrt{r}/2),
$$
we get 
\begin{equation}\label{e:plaques}
Kr^{d-1}\P(\ell \ge \sqrt{r}/2) \to 0  \mbox{ as } r \to \infty.
\end{equation}
Furthermore, 
\begin{align*}
 \P(H(r)) & \le \E(\card[\{(c,s) \in \xi : s \ge \sqrt{r} \mbox{ and } B(c,s) \cap B(0,2r) \neq \emptyset\}]) \\
& \le \lambda \int_{[\sqrt{r},\infty)} v_d (2r+s)^d \nu(ds) \\
& \le \lambda \int_{[\sqrt{r},\infty)} v_d (2s^2+s)^d \nu(ds) \\
\end{align*}
where $v_d$ denotes the volume of the unit ball of $\R^d$. Therefore 
\begin{equation} \label{e:recommencent}
\P(H(r)) \to 0 \mbox{ as } r \to \infty.
\end{equation}
The result follows from \eqref{e:poildechameau}, \eqref{e:plaques} and \eqref{e:recommencent}. \qed

\paragraph{Remark. } Replacing $\sqrt{r}$ by $ \alpha r/ \log r$ for some big enough $\alpha$ in the definition of the event $H(r)$, and replacing accordingly the event $\{ \ell \geq \sqrt{r} /2 \}$ by $\{ \ell \geq \log r / (2 \alpha) \}$, we could weakened the hypothesis $\E (R^{2d}) <\infty$ to $\E (R^d (\log R)^\beta) <\infty$ for some $\beta$, but only for the implication $2\Rightarrow 1$.

\paragraph{Proof of $1 \implies 3$.} This is a consequence of Theorem \ref{t:implicit}. \qed

\paragraph{Proof of $3 \implies 4$.} Note that $C$ is a subset of $B(0,D)$. As a consequence, $|C| \le v_d D^d$.
The result follows. \qed

\paragraph{Proof of $4 \iff 5$.} There exits actually simple inequalities between $\E(|C|)$ and $\E(\# C)$. 
Using Slivnyak’s theorem, see Proposition 4.1.1 in \cite{Moller}, we get
\begin{align*}
\E(\#C) 
 & = \E\left( \sum_{c \in \chi} 1_{c \in C} \right) \\
& = \lambda \int_{\R^d} dc \int_{(0,+\infty)} \nu(dr) \; \P\big(c \in C(\xi \cup (c,r))\big) 
\end{align*}
where $C(\xi \cup (c,r))$ denotes the connected component containing the origin in the Boolean model with the extra ball $B(x,r)$.
Thus,
\begin{align*}
\E(\#C) 
 & = \lambda \int_{\R^d} dc \int_{(0,+\infty)} \nu(dr) \; \P\big(B(c,r) \mbox{ touches } 0 \mbox{ or } B(c,r) \mbox{ touches } C\big) \\
 & = \lambda \int_{\R^d} dc \int_{(0,+\infty)} \nu(dr) \;  \P\big(c \in B(0,r) \cup (C + B(0,r)) \big) \\
 & = \lambda \E\big( |B(0,R) \cup (C+B(0,R))|\big)
\end{align*}
where $R$ is independent of $\xi$ and where the distribution of $R$ is $\nu$.
Thus,
\begin{equation} \label{e:retour}
\lambda  \E(C) \le \E(\#C).
\end{equation}
Moreover, using Lemma \ref{l:geo},
\begin{align}
 \E(\#C) & \le \lambda v_d \E(R^d) + \lambda \E\left(|C(1)+B(0,R)|\right)   \nonumber \\
  & \le \lambda v_d \E(R^d) + \lambda K^22^d \E([1+R]^d)\E(|C(1)+B(0,1)|). \label{e:lille}
\end{align}
But, by Lemma \ref{l:finalementfacile}, $\E(|C(1)+B(0,1)|)$ is finite. 
As $\E(R^{d})$ is also finite, the equivalence follows from \eqref{e:retour} and \eqref{e:lille}. \qed

\def\cprime{$'$} \def\cprime{$'$}


\begin{thebibliography}{10}

\bibitem{ATT2}
Daniel Ahlberg, Vincent Tassion, and Augusto Teixeira.
\newblock Sharpness of the phase transition for continuum percolation in $r^2$.
\newblock Available on Arxiv, 2016.

\bibitem{ATT}
Daniel Ahlberg, Vincent Tassion, and Augusto Teixeira.
\newblock Existence of an unbounded vacant set for subcritical continuum
  percolation.
\newblock Available on Arxiv, 2017.

\bibitem{Aizenman-Barsky}
Michael Aizenman and David~J. Barsky.
\newblock Sharpness of the phase transition in percolation models.
\newblock {\em Comm. Math. Phys.}, 108(3):489--526, 1987.

\bibitem{DC-video-17}
Hugo Duminil-Copin.
\newblock Sharp threshold phenomena in statistical physics.
\newblock Video available on the author's webpage, 2017.

\bibitem{DC-Raoufi-Tassion-17-2}
Hugo Duminil-Copin, Aran Raoufi, and Vincent Tassion.
\newblock Exponential decay of connection probabilities for subcritical voronoi
  percolation in $\mathbb{R}^d$.
\newblock Available on Arxiv, 2017.

\bibitem{DC-Raoufi-Tassion-17-1}
Hugo Duminil-Copin, Aran Raoufi, and Vincent Tassion.
\newblock Sharp phase transition for the random-cluster and potts models via
  decision trees.
\newblock Available on Arxiv, 2017.

\bibitem{DCT-sharp-Ising}
Hugo Duminil-Copin and Vincent Tassion.
\newblock A new proof of the sharpness of the phase transition for bernoulli
  percolation and the ising model.
\newblock {\em Communications in Mathematical Physics}, 343(2):725--745, 2016.

\bibitem{DCT-sharp-Bernoulli}
Hugo Duminil-Copin and Vincent Tassion.
\newblock A new proof of the sharpness of the phase transition for {B}ernoulli
  percolation on {$\Bbb{Z}^d$}.
\newblock {\em Enseign. Math.}, 62(1-2):199--206, 2016.

\bibitem{G-perco-boolean-model}
Jean-Baptiste Gou{\'e}r{\'e}.
\newblock Subcritical regimes in the {P}oisson {B}oolean model of continuum
  percolation.
\newblock {\em Ann. Probab.}, 36(4):1209--1220, 2008.

\bibitem{Gouere-Theret-17}
Jean-Baptiste Gouéré and Marie Théret.
\newblock Positivity of the time constant in a continuous model of first
  passage percolation.
\newblock {\em Electron. J. Probab.}, 22:21 pp., 2017.

\bibitem{last_penrose_2017}
Günter Last and Mathew Penrose.
\newblock {\em Lectures on the Poisson Process}.
\newblock Institute of Mathematical Statistics Textbooks. Cambridge University
  Press, 2017.

\bibitem{Meester-Roy-livre}
Ronald Meester and Rahul Roy.
\newblock {\em Continuum percolation}, volume 119 of {\em Cambridge Tracts in
  Mathematics}.
\newblock Cambridge University Press, Cambridge, 1996.

\bibitem{Menshikov-coincidence}
M.~V. Men{\cprime}shikov.
\newblock Coincidence of critical points in percolation problems.
\newblock {\em Dokl. Akad. Nauk SSSR}, 288(6):1308--1311, 1986.

\bibitem{Moller}
Jesper M{\o}ller.
\newblock {\em Lectures on random {V}orono\u\i\ tessellations}.
\newblock Springer-Verlag, New York, 1994.

\bibitem{Penrose17}
M.~D. {Penrose}.
\newblock {Non-triviality of the vacancy phase transition for the Boolean
  model}.
\newblock {\em ArXiv e-prints}, June 2017.

\bibitem{Schneider-Weil}
Rolf Schneider and Wolfgang Weil.
\newblock {\em Stochastic and integral geometry}.
\newblock Probability and its Applications (New York). Springer-Verlag, Berlin,
  2008.

\bibitem{vdB-BK}
J.~van~den Berg.
\newblock A note on disjoint-occurrence inequalities for marked {P}oisson point
  processes.
\newblock {\em J. Appl. Probab.}, 33(2):420--426, 1996.

\bibitem{Ziesche}
S.~{Ziesche}.
\newblock {Sharpness of the phase transition and lower bounds for the critical
  intensity in continuum percolation on $\mathbb{R}^d$}.
\newblock {\em ArXiv e-prints}, July 2016.

\end{thebibliography}
\end{document}